\documentclass[12pt]{amsart}
\def\date{7 July 2010}
\usepackage{verbatim}
\usepackage{amsmath}
\usepackage{amssymb}
\usepackage{amscd}
\usepackage{graphicx}
\usepackage{psfrag}

\numberwithin{equation}{section}

\newtheorem{thm}[equation]{Theorem}
\newtheorem{lemma}[equation]{Lemma}

\newtheorem*{theorem*}{Theorem}

\theoremstyle{definition}
\theoremstyle{remark}

\def\proof{{\noindent \bf Proof. }}

\def\qed{\hfill$\square$\bigskip\medskip}

\def\PP{{\mathbb P}}

\begin{document}

\title[minor-minimal planar graphs of even branch-width ]{ Minor-minimal planar graphs of even branch-width}
\author{%
  Torsten Inkmann
  and
  Robin Thomas
}
\address{School of Mathematics, Georgia Institute of Technology,
Atlanta, GA 30332-0160, USA}
\thanks{Partially supported by NSF grant DMS-0701077. \date.}

\begin{abstract}
Let $k\ge1$ be an integer, and
let $H$ be a  graph with no isolated vertices embedded in the projective plane 
such that every homotopically non-trivial closed curve intersects $H$
at least $k$ times,
and the deletion and contraction of any edge in this embedding
results in an embedding that no longer has this property.
%of face-width $k$,
Let $G$ be the planar double cover of $H$ obtained by lifting $G$
into the universal covering space of the projective plane, the sphere.
We prove that $G$ is minor-minimal of branch-width $2k$.
We also exhibit examples of  minor-minimal planar graphs of branch-width $6$
that do not arise this way.
\end{abstract}

\maketitle

\section{Introduction}
\label{SectionIntroduction}

\begin{figure}
\label{FigPetersenDC}
\centering
\includegraphics[totalheight=0.3\textheight]{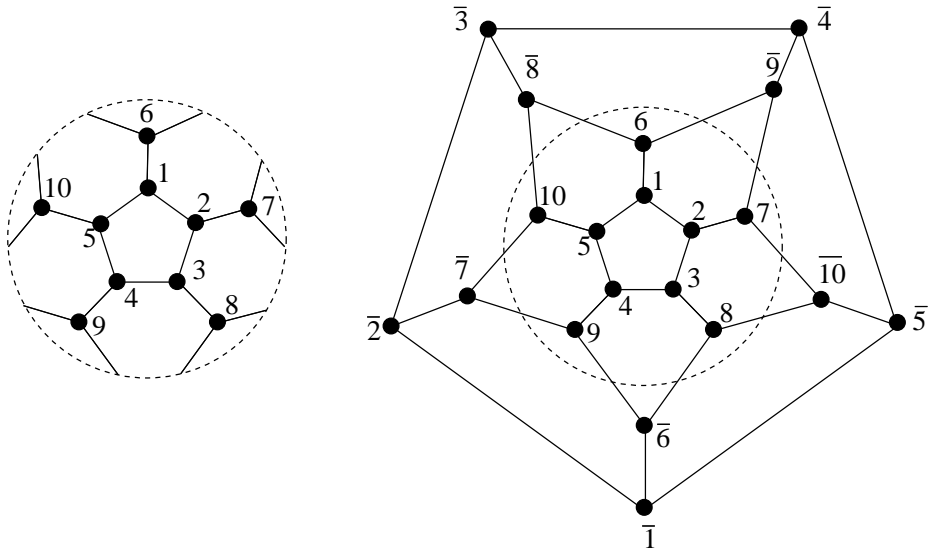}
\caption{An embedding of the Petersen graph in $\PP$ and its planar double cover, the Dodecahedron}
\end{figure}

The projective plane $\mathbb{P}$ is obtained from a closed disk by identifying
diagonally opposite pairs of points on the boundary of the disk. 
Given a graph $H$ embedded in $\mathbb{P}$ its {\em planar double cover}
$G$ is the lift of $H$ into the universal covering space of $\mathbb{P}$,
the sphere.
Thus to every vertex $v$ of $H$ there correspond two vertices $v_1,v_2$
of $G$; we say that $v$ is the {\em projection} of $v_1$ and $v_2$,
and that $v_1$ and $v_2$ are the {\em lifts} of $v$.
Similarly, we speak of projections and lifts of paths, cycles, walks,
and faces.
This construction is illustrated in Figure~\ref{FigPetersenDC},
where the Dodecahedron is shown to be a planar double cover of
the Petersen graph.
%(The planar double cover of the Petersen graph is actually unique, because 
%every two embeddings of the Petersen graph in the projective plane
%are related by an automorphism of the Petersen graph.)
In particular, if $W$ is a walk in $G$ with ends the two lifts of a vertex
$v\in V(H)$, then the projection of $W$ is a homotopically non-trivial
closed walk in $H$ with both ends $v$.

Let $H$ and $G$ be as in the above paragraph.
Our objective is to relate the representativity of $H$ (also known as face-width)
and the branch-width of 
$G$, two important parameters that we now review. 
A graph is a {\em minor} of another if the first can be obtained
from a subgraph of the second by contracting edges.
Let $k\ge0$ be an integer. 
We say that a graph $H$ embedded in $\mathbb{P}$ 
(a ``projective plane graph") is {\em $k$-representative}
if every homotopically
%The {\em face-width} of $H$, also known as {\em representativity},
% is the maximum integer $k$ such that every (homotopically) 
non-trivial closed curve in $\mathbb{P}$ intersects $H$ at least $k$ times.
(We mean all curves, including those passing through vertices
of $H$; in fact, it suffices to restrict oneself to curves that
intersect the graph only in vertices.)
This concept has received a lot of attention in the literature;
we refer to~\cite{MohTho} for more information. 
%We will need minor-minimal graphs of a given representativity.
We say that $H$ is {\em minor-minimal $k$-representative}
if $H$ is $k$-representative,
has no isolated vertices and for every edge $e$ of $H$
the embedded graphs obtained from $H$ by deleting and contracting $e$
are no longer $k$-representative.
(One could also define a related concept where we say that no
proper minor of $H$, taken as an abstract graph, has a $k$-representative
embedding in $\mathbb{P}$. For $k\ge3$ these two notions coincide,
but we use the former definition.)
Randby~\cite{Ran} proved that given two minor-minimal $k$-representative
graphs in $\mathbb P$,
each can be obtained from the other by means of
repeated application of $\Delta Y$- and $Y \Delta$-exchanges.
Since for each $k$ there is a natural example of a minor-minimal
$k$-representative projective plane graph, namely the $k\times k$
projective grid, Randby's result gives a convenient way to generate
all minor-minimal $k$-representative projective plane graphs.
Incidentally, Schrijver~\cite{Schmintorus} proved an analogue of Randby's result
for the torus, as well as a related result for
arbitrary orientable surfaces~\cite{Schkernels}.
%\nocite{SchHomotopic}\nocite{SchHomCirc}

A {\em branch-decomposition} of a graph $G$ is pair $(T, \eta)$, where 
$T$ is a tree with all vertices of degree one or three, and $\eta$
is a bijection between the leaves (vertices of degree one) of $T$ and $E(G)$. 
%A {\em separation} of $G$ is a pair of subgraphs $(A,B)$ of $G$ so that 
%$E(A) \cap E(B) = \emptyset$ and $A \cup B = G$. 
%The order of a separation $(A,B)$ is defined as $|V(A \cap B)|$. 
For $f\in E(T)$ let $T_1,T_2$ be the two components of $T\backslash f$, and
let $X_i$ be the set of leaves of $T$ that belong to $T_i$.
We define the {\em order} of $f$ to be the number of vertices of $G$
incident both with an edge in $\eta(X_1)$ and an edge in $\eta(X_2)$.
The {\em width} of   $(T, \eta)$ is the maximum order of an edge of $T$,
and the {\em branch-width} of $G$ is the minimum width of a
branch-decomposition of $G$, or 0 if $|E(G)|\le1$, in which case $G$
has no branch-decomposition.

Computing branch-width is NP-hard~\cite{SeyThoRat}, but there is
a polynomial time algorithm when $G$ is planar~\cite{SeyThoRat}. 
Thus one might expect that planar branch-width is better behaved
in other respects as well, such as in terms of excluded minors.
Since taking minors does not increase branch-width, as is easily seen,
graphs of branch-width at most $t$ are characterized by 
{\em excluded minors}, the list of minor-minimal graphs of branch-width
$t+1$. It follows from~\cite{RobSeyGM4} that this list is finite.

Our first result is that planar double covers of minor-minimal
$k$-representative projective plane graphs  are excluded minors for
odd branch-width:

\begin{thm}
\label{main thm}
Let $k \geq 1$ be an integer, 
let $H$ be a minor-minimal $k$-representative projective plane graph,
and let $G$ be the double cover of $H$.
Then $G$ is a minor-minimal graph of branch-width $2k$. 
\end{thm}

Thus, in view of Randby's theorem mentioned above, this gives
a description of a class of minor-minimal graphs of branch-width $2k$.
Since for $k=1$ and $k=2$ this class actually includes all minor-minimal
graphs of branch-width $2k$~\cite{BodThi, DhaPhD, JohRob},
\nocite{BodThi}
\nocite{DhaPhD}
\nocite{JohRob}
one might wonder whether this holds in general.
Unfortunately, it does not. 
In Section~\ref{SectionNotDoubleCover} we exhibit examples of minor-minimal
planar graphs of branch-width $6$ that are not double covers
of any graph.

The paper is organized as follows.
We prove the lower bound for Theorem~\ref{main thm} in Section~\ref{sec:lower}
and the corresponding upper bound in Section~\ref{sec:upper}.
The examples are presented in Section~\ref{SectionNotDoubleCover}.

The second author would like to acknowledge helpful conversations with 
P.~D.~Seymour from the summer of 1990; in particular, the question
whether all minor-minimal graphs of even branch-width arise as in
Theorem~\ref{main thm} was inspired by those conversations.
The second author would also like to acknowledge that the results of this paper
appeared (in a slightly different form) in the PhD 
dissertation~\cite{InkPhD} of the first author.
%We are indebted to Ilya Hicks for letting us use his implementation
%of the algorithm from~\cite{SeyThoRat}.

%%%%%%%%%%%%%%%%%%%%%%%%%%%%%%%%%%%%%%%%%%%%%%%%%%%%%%%%%%%%%%%%%%%%%%%%%%%%

\section{Lower Bound}
\label{sec:lower}
We will make use of the result of~\cite{SeyThoRat} that the branch-width
of a planar graph $G$ is equal to one half the ``carving-width" of
the ``medial graph" of $G$. Here are the definitions.
Let $G$ be a graph drawn in a surface in such a way that every face of $G$ is
homeomorphic to an open disk.  Let us choose, for each vertex
$v\in V(G)$, one of the two cyclic orderings of edges incident
with $v$ (with each loop occurring twice in the ordering) and
designate it as {\em clockwise}.  By an {\em angle} at $v$ we mean
an ordered pair $(e,e')$ of edges incident with $v$ such that $e'$ 
immediately follows $e$ in the clockwise order around $v$.  To each 
angle there naturally corresponds a face $f$ of $G$ incident with $e$ and $e'$.
The {\em medial graph} of $G$ is the graph $M$ defined as follows.
For each edge $e\in E(G)$ choose a vertex $x_e$ positioned in the interior 
of $e$, and for each angle $(e,e')$ choose an edge joining $x_e$ and $x_{e'}$
inside the face of $G$ that corresponds to the angle $(e,e')$ in a small
neighborhood of $e$ and $e'$.  There is a certain ambiguity when the same
face corresponds to different angles, but there is a natural interpretation
under which the medial graph is unique up to homotopic shifts of edges.
Each face of $M$ corresponds to either a unique vertex of $G$ or a
unique face of $G$, and if $G$ has no loops or cut edges, then every
face of $M$ is bounded by a cycle.

A {\em carving} in a graph $G$ is a pair $(T, \eta)$, where $T$ is a tree
with all vertices of degree one or three, and
$\eta$ is a bijection from the leaves of $T$ to $V(G)$. 
Similarly as in a branch-decomposition, each edge $f$ of $T$ determines a
cut in $G$.
The {\em width} of  $(T, \eta)$ is the maximum order of those cuts,
over all $f\in E(T)$.
The {\em carving-width} of $G$ is the minimum width of a carving in $G$.
%, denoted by $cw(G)$, is defined to be the minimum width over all 
%carvings in $G$.
The following is~\cite[Theorem~7.2]{SeyThoRat}.

\begin{thm}
\label{bwcw}
Let $G$ be a connected plane graph with at least two edges,
and let $M$ be its medial graph. Then the branch-width of $G$
is half the carving-width of $M$.
\end{thm}

Thus it suffices to bound the carving-width of the medial graph of $G$
from Theorem~\ref{main thm}.
We will use the concept of an ``antipodality," introduced in~\cite{SeyThoRat}.
Let $G$ be a connected plane graph with planar dual $G^*$,
let $F(G)$ be the faces of $G$,
and let $k\ge0$ be an integer.
An {\em antipodality} in $G$ of range $\geq k$ is a function $\alpha$ with 
domain $E(G) \cup F(G)$, such that for all $e \in E(G)$, $\alpha(e)$ is a 
subgraph of $G$ with at least one vertex
and for all $f \in F(G)$, $\alpha(f)$ is a 
non-empty subset of $V(G)$, satisfying:
\begin{itemize}
\item[(A1)] If $e \in E(G)$, then no end of $e$ belongs to $V(\alpha(e))$
\item[(A2)] If $e \in E(G)$, $f \in F(G)$, and $e$ is incident with $f$, then $\alpha(f) \subseteq V(\alpha(e))$, and every component of $\alpha(e)$ has a vertex in $\alpha(f)$
\item[(A3)] If $e_1 \in E(G)$ and $e_2 \in E(\alpha(e_1))$ then every closed walk of $G^*$ using $e^*_1$ and $e^*_2$ has length $\geq k$.
\end{itemize}

The following is a special case of~\cite[Theorem~4.1]{SeyThoRat};
we will only need the easier ``if" part.

\begin{thm}
\label{thm:rat}
Let $M$ be a connected plane graph on at least two vertices,
and let $k \geq 0$ be an integer. 
Then $M$ has carving-width at least $k$ if and only if either 
some vertex of $M$ has degree at least $k$,
or $M$ has an antipodality of range $\geq k$.
\end{thm}

As we will see at the end of the next section, the following lemma
and  Theorems~\ref{bwcw} and~\ref{thm:rat} imply that
the planar double cover of
a minor-minimal $k$-representative projective plane graph
has branch-width at least $2k$.

%We are now ready to prove that a planar double cover of 
%a minor-minimal $k$-representative projective planar graph
%has branch-width at least $2k$.
%The result follows from Theorems~\ref{bwcw} and~\ref{thm:rat},
%and the following lemma.

\begin{lemma}
\label{lem:antipodality}
Let $k\ge2$ be an integer, 
let $H$ be a  $k$-representative projective plane graph, let
$G$ be a planar double cover of $H$, and let $M$ be the medial graph of $G$. 
Then $M$ has an antipodality of range $\ge 4k$. 
\end{lemma}

\proof
First we notice that since $k\ge2$ 
we may assume (by considering a subgraph of $H$) that
the graph $G$ has no loops or
cut edges.
For $v\in V(G)$ and a face $f\in F(G)$ incident with $v$ there is
a unique edge $a_{v,f}$ of $M$ incident with $v^*$ and $f^*$,
where $v^*$ and $f^*$ denote the 
faces of $M$ corresponding to $v$ and $f$, respectively.
For $v\in V(G)$ let $v'\in V(G)$ be the other vertex of $G$ with the
same projection as $v$, and for $f\in F(G)$ let $f'$ be defined
analogously.
Finally, for $x\in V(G)\cup F(G)$ let $C_x$ denote the cycle
bounding the face of $M$ corresponding to $x$.
We note that  $C_x$ is indeed a cycle, because $G$ has no loops or
cut edges.
We define the function $\alpha$ as follows: 
\begin{align*}
\alpha(a_{v,f}) &=  C_{{v'}} \cup C_{{f'}} 
\text{ for every } v\in V(G) \text{ and } f\in F(G)\\
\alpha(f) &= V(C_{{f'}}) \text{ for every } f \in F(M)
\end{align*}

We claim that $\alpha$ defines an antipodality of range $4k$ in $M$.
Conditions (A1) and (A2) follow easily from the fact that $G$ has no loops
or cut edges and that each $\alpha(a_{v,f})$ is connected.

To prove (A3) suppose for a contradiction that there is a closed
walk $W$ in $M^*$  of length at most $4k-1$
that uses $a^*$ and $b^*$ for some
$a=a_{v,f}\in E(M)$ and $b\in E(\alpha(a_{v,f}))$.
Thus $W$ uses both $v^*$ and $f^*$, and at least one of $(v')^*$, $(f')^*$.
In either case, $W$ includes a subwalk $W'$ of length at most $2k-1$ from
$x^*$ to $(x')^*$ for some $x\in V(G)\cup F(G)$.
Now the projection of $W'$ is a homotopically non-trivial closed walk
of length at most $2k-1$ in the dual of the medial graph of $H$.
That, in turn, implies that $H$ is not $k$-representative,
a contradiction. Hence $\alpha$ is an antipodality, as desired.
\qed

%%%%%%%%%%%%%%%%%%%%%%%%%%%%%%%%%%%%%%%%%%%%%%%%%%%%%%%%%%%%%%%%%%%%%%%%%%%%
\section{Upper bound}
\label{sec:upper}

In this section we show that if $G$ is as in Theorem~\ref{main thm},
then every proper minor of $G$ has branch-width at most $2k-1$.
Again, we find it convenient to work with carving-width of the
medial graph.
Thus we need to translate minimal $k$-representativity into the language
of medial graphs.
Let $G$ be a graph drawn in a surface, 
let  $v$ be a vertex of $G$ of degree four, let $e_1,e_2,e_3,e_4$ be the
four edges of $G$ incident with $v$ listed in the cyclic order of appearance
around $v$, let $v_i$ be the other end of $e_i$,
and let $f_i$ be the face of $G$ incident with
$e_i$ and $e_{i-1}$, where $e_0$ means $e_4$.
Let $G'$ be obtained from $G\backslash v$ by inserting two edges,
one with ends $v_1$ and $v_4$ in the face $f_1$ and the other with
ends $v_2$ and $v_3$ in the face $f_3$.
We say that $G'$ was obtained from $G$ by 
{\em opening at $v$ through the faces $f_2$ and $f_4$}, or simply by
{\em opening at $v$}.
Let $k\ge 1$ be an integer.
A $4$-regular graph $G$ drawn in $\mathbb P$ is {\em $k$-tight} if
%there exists an integer $k\ge 1$ such that 
every homotopically
non-trivial closed walk in $G^*$ has length at least $k$, and for
every graph $J$ obtained from $G$ by opening at some vertex there
exists a homotopically
non-trivial closed walk in $J^*$ of length at most $k-1$.
We say that $G$ is {\em tight} if it is $k$-tight for some 
integer $k\ge1$.
The following is shown in~\cite{Schkernels} and is also easy to see.

\begin{lemma}
\label{fwtight}
Let $k\ge1$ be an integer, let $G$ be a connected graph in the projective plane,
and let $M$ be the medial graph of $G$.
Then $G$ is a minor-minimal $k$-representative graph in $\mathbb P$ 
if and only if $M$ is $2k$-tight.
\end{lemma}

We will need the following characterization of tight graphs
in terms of straight ahead decompositions.
Let $G$ be a graph with all vertices of degree four or one drawn in a surface.
Let $v$ be a vertex of $G$ of degree four, and let $e_1,e_2,e_3,e_4$
be the edges of $G$ incident with $v$ listed in the order in which
they appear around $v$. We say that the edges $e_1$ and $e_3$ are
{\em opposite}. 
Let $F_1,F_2,\ldots,F_r$ be the equivalence classes of the
transitive closure of the opposite relation, and let
$G_i$ be the subgraph of $G$ with edge-set $F_i$ and vertex-set
all vertices incident with edges of $F_i$.
We say that $G_1,G_2,\ldots,G_r$ is the {\em straight ahead decomposition}
of $G$.
The next result follows from a theorem of Lins~\cite{Linpp}.
Schrijver~\cite{SchHomotopic} obtained an analogous result for
orientable surfaces.

\begin{thm}
\label{lins}
Let $G$ be a $4$-regular graph drawn in the projective plane.
Then $G$ is tight if and only if the straight ahead decomposition
of $G$ consists of homotopically non-trivial cycles such that
every two intersect exactly once.
\end{thm}

As usual, a walk in a graph $G$ is an alternating sequence of 
vertices and edges of $G$. It has an {\em origin} and {\em terminus},
called its {\em ends}. A walk is {\em closed} if its ends are equal.
Even though any vertex of a closed
walk may be regarded as its origin and terminus,
for our purposes changing the ends results in a different walk.
This subtle point will be important later.

Let $G$ be a graph drawn in a surface, let $f,f'$ be two faces
of $G$, and let $W_1,W_2$ be two walks in $G^*$ with origin $f^*$
and terminus $(f')^*$.
Let $v$ be a vertex of $G$ of degree four, let $e_1,e_2,e_3,e_4$ be the 
four edges of $G$ incident with $v$ listed in the cyclic order of appearance
around $v$, and let $f_i$ be the face of $G$ incident with 
$e_i$ and $e_{i-1}$, where $e_0$ means $e_4$.
If $W_1$ includes the subwalk $f_1^*,e_1^*,f_2^*,e_2^*,f_3^*$,
and $W_2$ is obtained from $W_1$ by replacing (one occurence of)
that subwalk by the walk
$f_1^*,e_4^*,f_4^*,e_3^*,f_3^*$,
then we say that $W_2$ was obtained from $W_1$ by a 
{\em $\Delta\nabla$-exchange}.
We write $W_1=W_2*v$; then also $W_2=W_1*v$.
Let us remark the obvious fact that $W_1$ and $W_2$ have the
same origin and terminus.

Let $\Delta$ denote the closed unit disk in $\mathbb{R}^2$.
Let $G$ be a graph drawn in $\Delta$ such that
the vertices $v_1,v_2,\ldots,v_{2k}$ and only these vertices
are drawn on the boundary of $\Delta$ in the order listed,
where $v_1,v_2,\ldots,v_{2k}$ have degree one and all
other vertices of $G$ have degree four.
Assume further that the straight ahead decomposition of $G$ is of
the form
$P_1,P_2,\ldots,P_k$, where each $P_i$ is a path with one end
in $\{v_1,v_2,\ldots,v_k\}$ and the other end in 
$\{v_{k+1},v_{k+2},\ldots,v_{2k}\}$.
Finally, assume that for $i\ne j$ the paths $P_i$ and $P_j$ intersect
in at most one vertex.
In those circumstances we say that $G$ is a {\em graft}, and that
$v_1,v_2,\ldots,v_{2k}$ are its {\em attachments}.
For $i=1,2,\ldots,2k$ let 
$e_i$ denote the unique edge of $G$ incident with $v_i$, and let
$f_i$ be the face of $G$ incident with
$e_{i-1}$ and $e_i$, where $e_0$ means $e_{2k}$.
Let us emphasize that $G$ is embedded in $\Delta$, and hence 
$f_i$ and $f_{i+1}$ are distinct faces, each incident with a segment
of the boundary of $\Delta$.
Let $G^*$ denote the geometric dual of $G$.
By a {\em broom} in $G$ we mean a walk in $G^*$ 
from $f_1^*$ to $f_{k+1}^*$
of length exactly $k$.
(The existence of $P_1,P_2,\ldots,P_{k}$ implies that every
walk in $G^*$
from $f_1^*$ to $f_{k+1}^*$ has length at least $k$.)
There are two natural examples of brooms, namely the walks with
edge-sets $f_1^*,f_2^*,\ldots,f_{k+1}^*$ and
$f_{1}^*,f_{2k}^*,f_{2k-1}^*,\ldots,f_{k+1}^*$.
Those two brooms will be called the {\em extreme brooms} of $G$.

Let $G$ be a graph drawn in a surface with every vertex of
degree four or one, let $f,f'$ be two
faces of $G$, and let $W_0$ and $W$ be two walks in $G^*$ with origin
$f$ and terminus $f'$.
We say that $G$ is {\em sweepable} from $W_0$ to $W$ if the
vertices of $G$ of degree four can be numbered $v_1,v_2,\ldots,v_n$
so that $W_i:=W_{i-1}*v_i$ is well-defined for all $i=1,2,\ldots,n$
and $W_n=W$.

\begin{lemma}
\label{graft}
Let $k\ge1$ be an integer,
let $G$ be a graft with attachments $v_1,v_2,\ldots,v_{2k}$,
and let $W,W'$ be the extreme brooms of $G$.
Then $G$ is sweepable from $W$ to $W'$.
\end{lemma}

\proof
Let $P_1,P_2,\ldots,P_k$ be the straight ahead decomposition of $G$,
numbered so that $v_i$ is an end of $P_i$ for $i=1,2,\ldots,k$.
If the paths $P_1,P_2,\ldots,P_k$ are pairwise disjoint, then
$G$ is a matching, and hence
$W=W'$ and the theorem holds. Thus we may assume that some two
of the paths $P_i$ intersect, and we proceed by induction on $|V(G)|$.
By a {\em wedge} we mean an ordered pair $(i,j)$ of distinct integers
from $\{1,2,\ldots,k\}$ such that the paths $P_i$ and $P_j$
intersect in a (unique) vertex $v$, and $v_iP_iv$, the subpath
of $P_i$ with ends $v_i$ and $v$, is not intersected by any other
$P_{j'}$. Since some two paths $P_i$ intersect, there exists a wedge,
and so we may select a wedge $(i,j)$ with $|i-j|$ minimum.

We claim that $|i-j|=1$. To prove this claim we may assume that
$i+1<j$. By planarity the path $P_{i+1}$ intersects $P_i\cup P_j$,
and hence there exists a wedge $(i+1,j')$. We deduce from the minimality 
of $|i-j|$ that $j'\not\in\{i,i+1,\ldots,j\}$, and hence $P_{j'}$
intersects $v_iP_iv\cup vP_jv_j$ by planarity. In fact, it intersects
$v_iP_iv\cup vP_jv_j$ an even number of times.
But $P_{j'}$ does not intersect  $v_iP_iv$, because $(i,j)$ is a wedge,
% by the definition of a wedge,
and hence  $P_{j'}$ intersects $P_j$ at least twice, contrary to the definiton
of a graft.
This contradiction proves our claim that $|i-j|=1$.
It follows that both $v_iP_iv$ and $v_jP_jv$ have length one.

Thus we can apply a $\Delta\nabla$-exchange at the vertex $v$
to one of the brooms $W,W'$,
say to $W'$, to obtain a broom $W'':=W'*v$.
Let $x_i$ be the neighbor of $v$ on $P_i$ other than $v_i$, and let
$x_j$ be defined analogously. 
Let $G'$ be obtained from $G$ by deleting $v,v_i,v_j$ and adding two
new vertices of degree one, joined to $x_i$ and $x_j$, respectively.
Then $G'$ with its natural drawing in a disk forms a graft, and $W$
and $W''$ can be regarded as the two extreme brooms in $G'$.
By the induction hypothesis the graft $G'$ is sweepable from $W$
to $W''$; by considering the corresponding ordering of vertices of
$G'$ and appending $v$ at the end we obtain a desired ordering 
of the vertices of $G$, showing that $G$ is sweepable from $W$ to $W'$.
\qed

\begin{thm}
\label{thm:pp}
Let $k\ge1$ be an integer,
let $H$ be a $4$-regular $k$-tight graph in $\mathbb{P}$,
%let $H$ be a graph drawn in $\mathbb{P}$ such that its
%straight ahead decomposition consists of  $k$ homotopically
%non-trivial cycles any two of which intersect exactly once,
and let $W$ be a homotopically non-trivial closed walk in $H^*$
of length $k$.
Then $H$ is sweepable from $W$ to $W$.
\end{thm}

\proof
By Theorem~\ref{lins} the straight ahead decomposition of $H$
consists of $k$ homotopically non-trivial cycles such that every
two of them intersect exactly once.
We cut $H$ open along $W$ and construct a graft $G$ as follows. 
For every edge $e^*$ of $W$ we cut the corresponding edge $e$ into
two by inserting two new vertices of degree one ``in the middle of" $e$.
The theorem follows by applying Theorem~\ref{graft} to the
resulting graft $G$.
\qed

\begin{lemma}
\label{planarcover}
Let $k\ge1$ be an integer,
let $H$ be a $4$-regular $k$-tight graph in $\mathbb{P}$,
%let $H$ be a graph drawn in $\mathbb{P}$ such that its
%straight ahead decomposition consists of  $k$ homotopically
%non-trivial cycles any two of which intersect exactly once,
and let $G$ be the planar double cover of $H$.
Let $f_1,f_2$ be the two lifts of some face $f$ of $H$, and
let $W_1$ be a walk of length $k$ in $G^*$ with origin $f_1$ and 
terminus $f_2$.
Then $G$ is sweepable from $W_1$ to $W_1$.
\end{lemma}

\proof
Let $W$ be the projection of $W_1$, and let $W_2$ be the other lift of $W$.
Then $W$ is a homotopically non-trivial closed walk in $H^*$.
By Theorem~\ref{thm:pp} $H$ is sweepable from $W$ to $W$;
let $v_1,v_2,\ldots,v_n$ be the corresponding ordering of the
vertices of $H$.
Then  each $v_i$ has a lift $v_i^1\in V(G)$ such that 
$W_2=W_1*v_1^1*v_2^1*\cdots*v_n^1$.
Now let $v_i^2$ be the other lift of $v_i$.
It follows that 
$W_1=W_2*v_1^2*v_2^2*\cdots*v_n^2$,
and hence the sequence $v_1^1,v_2^1,\ldots,v_n^1,v_1^2,v_2^2,\ldots,v_n^2$
shows that $G$ is sweepable from $W_1$ to $W_2$.
\qed

\begin{lemma}
\label{open}
Let $k\ge1$ be an integer,
let $H$ be a $4$-regular $k$-tight graph in $\mathbb{P}$,
%let $H$ be a graph drawn in $\mathbb{P}$ such that its
%straight ahead decomposition consists of  $k$ homotopically
%non-trivial cycles any two of which intersect exactly once, 
and let $G$ be the planar double cover of $H$.
Let $G_1$ be obtained from $G$ by opening at $\overline u\in V(G)$.
Then $G_1$ has carving-width at most $2k-1$.
\end{lemma}

\proof
Let $\overline e_1,\overline e_2 ,\overline e_3 ,\overline e_4 $
%Let $\overline{e_1},\overline{e_2},\overline{e_3},\overline{e_4}$
be the edges of $G$ incident with $\overline u$ listed in cyclic
order around $\overline u$, and let $\overline{f_i}$ be the
face of $G$ incident with $\overline e_i $ and $\overline e_{i-1} $,
where $\overline e_0 $ means $\overline e_4 $.
Let the numbering be such that $G_1$ is obtained by opening at
$\overline u$ through the faces $\overline{f_2}$ and $\overline{f_4}$.
Let $u, e_i,f_i$ be the projections of 
$\overline u, \overline e_i, \overline{f_i}$, respectively.
Let $H_1$ be the graph obtained from $H$ by opening at $u$
through the faces $f_2$ and $f_4$, and let $f_{24}$ be the resulting
new face.
Since $H$ is $k$-tight there exists a homotopically non-trivial closed
walk $Z_1$ in $H_1^*$ of length $k-1$ with origin and terminus
$f_{24}^*$. Since $H$ is $k$-tight, we deduce that $f^*_{24}$ is the only
repeated vertex in $Z_1$.
To $Z_1$ there corresponds a closed walk $Z$ in $H^*$ starting with
$f_2^*e_2^*f_3^*e_3^*f_4^*$ and ending in $f_2^*$.
Let $\overline Z$ be the lift of $Z$ that starts with 
$(\overline{f_2})^*(\overline e_2)^*(\overline{f_3})^*
(\overline e_3)^*(\overline{f_4})^*$ 
and ends in $(\tilde{f_2})^*$, where $\tilde f_2$ is the other lift of $f_2$.
There is a corresponding walk $\overline Z_1$
in $G_1^*$ with origin $(\overline f_{24})^*$
and terminus $(\tilde{f_2})^*$, where $\overline f_{24}$
is the face resulting from the opening of $G$ at $\overline u$
that creates $G_1$.

By Lemma~\ref{planarcover} the graph $G$ is sweepable
from $\overline Z$ to $\overline Z$; let $v_1,v_2,\ldots,v_n$
be the corresponding ordering of the vertices of $G$, and let
$W_0:=\overline Z$ and $W_i:=W_{i-1}*v_i$ be the corresponding walks.
Since each $W_i$ starts at $(\overline{f_2})^*$ we deduce that the
$\Delta\nabla$-exchange at $\overline u$ replaces
$(\overline{f_2})^*(\overline e_1)^*(\overline{f_1})^*
(\overline e_4)^*(\overline{f_4})^*$
by
$(\overline{f_2})^*(\overline e_2)^*(\overline{f_3})^*
(\overline e_3)^*(\overline{f_4})^*$
or vice versa. 
From the symmetry we may assume that it replaces the former by the latter;
then we may further assume that $\overline u=v_n$.
We deduce that for all $i=1,2,\dots, n-2$ the concatenation of $\overline Z_1$ and $W_i$
separates $\{v_1,v_2,\ldots,v_{i}\}$ from $\{v_{i+1},v_{i+2},\ldots,v_{n-1}\}$
in $G_1$.
Now let $T$ be the tree obtained from a path with vertices
$r_1,r_2,\ldots,r_{n-1}$ in order by adding, for
each $i=2,3,\ldots,n-2$, a new vertex $t_i$ and joining it by
an edge to $r_i$. Let $t_1=r_1$ and $t_{n-1}=r_{n-1}$, and
let $\eta(t_i)=v_i$.
Since $\overline Z_1$ has length $k-1$ and each $W_i$ has length $k$
we deduce that $(T,\eta)$ is a carving decomposition of $G_1$ of 
width at most $2k-1$, as desired.~\qed

%We are now ready to prove Theorem~\ref{main thm}.

\noindent
{\bf Proof of Theorem~\ref{main thm}}.
Let $k\ge1$ be an integer, let $H$ be a minor-minimal $k$-representative
graph in $\mathbb P$, and let $G$ be the double cover of $H$.
If $k=1$, then $H$ consists of one vertex and one edge, and hence
$G$ has two vertices and two edges between them. It follows that
the theorem holds in that case, and hence we may assume that $k\ge 2$.

We first show that $G$ has branch width at least $2k$.
Let $M$ be the medial graph of $G$.
By Lemma~\ref{lem:antipodality} the graph $M$ has an antipodality of
range $\ge 4k$, and hence has carving-width at least $4k$ by Theorem~\ref{thm:rat}.
It follows from Theorem~\ref{bwcw} that $G$ has branch-width at least $2k$,
as desired.

Let $G_1$ be obtained from $G$ by deleting or contracting an edge $e$.
To complete the proof we show that $G_1$ has branch-width at most $2k-1$;
that will imply that the branch-width of $G$ is exactly $2k$ and that
it is minor-minimal.
To this end let $M_1$ be the medial graph of $G_1$; then $M_1$ is
obtained from $M$ by opening at the vertex of $M$ that correponds to $e$.
Let $N$ be the medial graph of $H$; 
then $N$ is $2k$-tight by Lemma~\ref{fwtight} and $M$ is the planar
double cover of $N$.
By Lemma~\ref{open} the graph $M_1$ has carving-width at most $4k-2$,
and hence $G_1$ has branch-width at most $2k-1$ by Theorem~\ref{bwcw},
as desired.~\qed

%%%%%%%%%%%%%%%%%%%%%%%%%%%%%%%%%%%%%%%%%%%%%%%%%%%%%%%%%%%%%%%%%%%%%%%%%%%
\section{Minor-minimal planar graphs which are not double covers}
\label{SectionNotDoubleCover}

\begin{figure}[ht]
\psfrag{n=15}{}
\psfrag{n=14}{}
\psfrag{m=33}{}
\psfrag{m=32}{}
\psfrag{m=31}{}
\psfrag{r=18}{}
\psfrag{r=19}{}
\psfrag{r=20}{}
\centering
\includegraphics[totalheight=0.5\textheight]{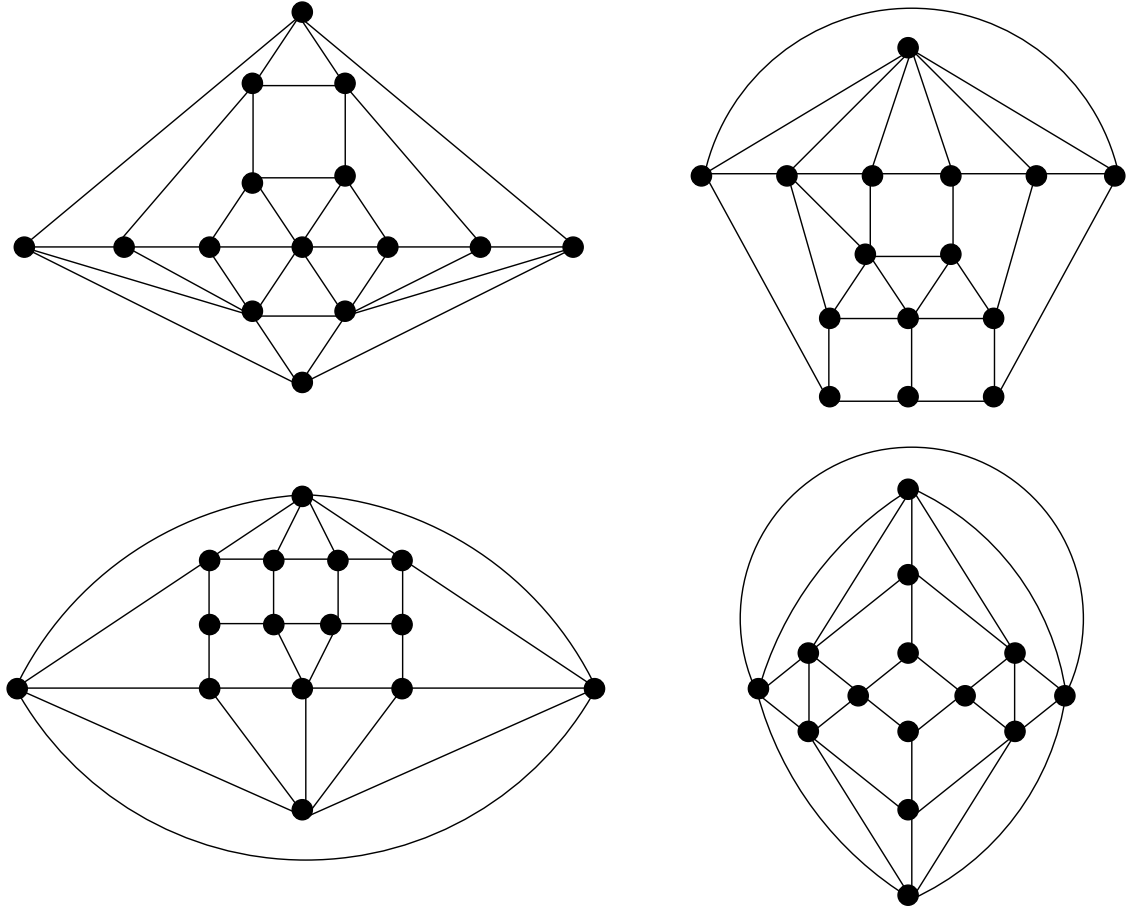}
\caption{Minor-minimal graphs of branch-width 6 which are not double covers}
\label{fig:FigBw6min}
\end{figure}

In~\cite{Barpoly} and~\cite{Vitpp} the authors found all the seven minor-minimal
$3$-representative embeddings of graphs in the projective plane.
(There is a possible confusion here: these are embeddings of six different
graphs, because one of those graphs has two different 
$3$-representative embeddings in the projective plane. 
The six graphs are all members of the so-called Petersen family.
The seventh member of the family, the graph obtained from $K_{4,4}$ by
deleting one edge, does not embed in the projective plane.)
By Theorem~\ref{main thm} the seven minor-minimal $3$-representative
projective planar embeddings give rise to
seven non-isomorphic minor-minimal planar graphs of branch-width $6$.
%$6$ which arise as planar double covers of the seven minor-minimal 
%$3$-representative embeddings of graphs in the projective plane 
%(the latter were determined in 
%\cite{Barpoly} and \cite{Vitpp}).
However, Figure~\ref{fig:FigBw6min} shows four additional
minor-minimal planar graphs of branch-width $6$, and their geometric
duals provide four additional examples.
In fact, none of these graphs is a double cover of any graph,
because none of them has a fixed point-free automorphism of order two.
We have generated those graphs using a computer program,
which used as a subroutine an implementation of the algorithm 
from~\cite{SeyThoRat} written by Hicks.
We are indebted to Ilya Hicks for letting us use his program.

Let $k\ge1$ be an integer, and let $H,G$ be as in Theorem~\ref{main thm}.
Let $v\in V(H)$ have degree three, and let $H'$ be obtained from $H$
by a $Y\Delta$-exchange at $v$; that is, $H'$ is obtained from $H$
by deleting $v$ and adding an edge joining every pair of neighbors of $v$.
The graph $H'$ gives rise to a minor-minimal graph $G'$ of branch-width $2k$
as in  Theorem~\ref{main thm}, where $G'$ can be regarded as having been
obtained from $G$ by a $Y\Delta$-exchange at both the lifts of $v$,
say $v_1$ and $v_2$.
Now a question arises about the status of the graph $G_1$ obtained from
$G$ by performing a $Y\Delta$-exchange at $v_1$ only.
For instance, when $k=2$ the graph $G$ could be the cube, but in that case
the graph $G_1$ has branch-width $3$.
We suspect that in general the graph $G_1$ has branch-width $2k-1$,
and that it can be shown using similar method as our proof of 
Theorem~\ref{main thm}.
However, it does not seem to follow easily from anything we have done,
and so we do not pursue this question further.

\begin{figure}[ht]
\centering
\includegraphics[totalheight=0.3\textheight]{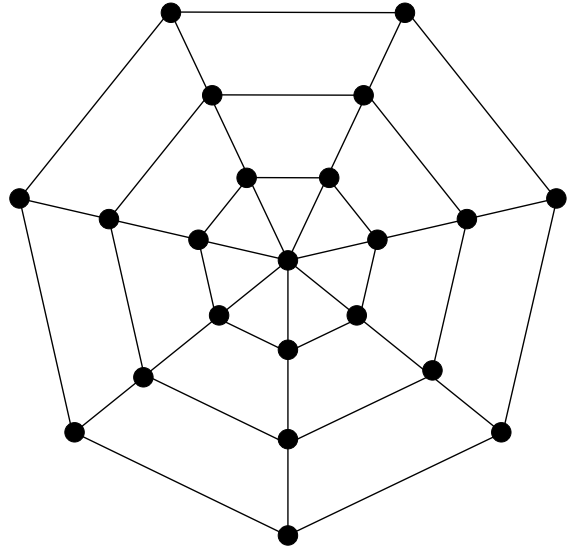}
\caption{The graph $G_{3}$}
\label{fig:FigCobweb3}
\end{figure}

Finally, here is an example of a minor-minimal planar graph of odd branch-width.
For every $k \geq 2$ consider the planar $k \times (2k+1)$ circular grid 
with $k$ concentric cycles and $2k+1$ paths joining the cycles, 
and add a new vertex that is 
connected to all $2k+1$ vertices on the innermost cycle to obtain a graph 
$G_k$. The graph $G_3$ is depicted in Figure~\ref{fig:FigCobweb3}.
It can be shown that  $G_{k}$ has branch-width exactly $2k+1$ and is 
minor-minimal with that property.

%\section{Conclusion}
%\label{SectionConclusion}
%method can be turned into polytime algo%
%can other minormin planar graphs of bw=2k be characterized?%
%what about bw=2k+1?%

\bibliographystyle{plain}
%\bibliography{biblio}

\begin{thebibliography}{99}

%\input refs.tex


\bibitem{Barpoly}D. Barnette, The minimal projective plane polyhedral maps,
{\it Applied geometry and discrete mathematics,}  63--70, 
DIMACS Ser. Discrete Math. Theoret. Comput. Sci., 4, 
Amer. Math. Soc., Providence, RI, 1991.

\bibitem{BodThi}
H.~L. Bodlaender and D.~M. Thilikos,
Graphs with branchwidth at most three,
\newblock {\em J. Algorithms}, {\bf32} (1999), 167--194.

\bibitem{DhaPhD}
J.~Dharmatilake,
Binary Matroids With Branch-Width Three,
PhD thesis, Ohio State University, 1994.

\bibitem{InkPhD} T.~Inkmann,
Tree-based decompositions of graphs on surfaces and applications 
to the Traveling Salesman Problem,
Ph.D.\ Dissertation, Georgia Institute of Technology, April 2008.

\bibitem{JohRob}
E.~Johnson and N.~Robertson,
private communication.

\bibitem{Linpp} S.~Lins,
A minimax theorem on circuits in projective graphs,
{\it\JCTB} {\bf 30} (1981), 253--262.

\bibitem{MohTho} B.~Mohar and C.~Thomassen, Graphs on surfaces, 
Johns Hopkins University Press, Baltimore, MD, 2001.

\bibitem{Ran} S.~P.~Randby,
Minimal embeddings in the projective plane,
{\it\JGT} {\bf25} (1997), 153--164.

\bibitem{RobSeyGM4}N.~Robertson and P.~D.~Seymour, 
Graph Minors IV. Tree-width and well-quasi-ordering,  
\JCTB\ {\bf 48} (1990), 227--254

\bibitem{SchHomotopic} A. Schrijver, 
Decomposition of graphs on surfaces and a homotopic circulation theorem, 
{\it\JCTB} {\bf51} (1991), 161--210.

\bibitem{Schkernels} A. Schrijver, 
On the uniqueness of kernels, 
{\it\JCTB} {\bf55} (1992), 146--160.

\bibitem{Schmintorus} A.~Schrijver, 
Classification of minimal graphs of given representativity on the torus, 
{\it\JCTB} {\bf 61} (1994), 217--236.

\bibitem{SeyThoRat}P.~D.~Seymour and R.~Thomas, 
Call routing and the ratcatcher,
\COM\ {\bf14} (1994), 217--241.

\bibitem{Vitpp} R.~P.~Vitray,
Representativity and flexibility on the projective plane,
{\it Graph structure theory (Seattle, WA, 1991)},  341--347, 
Contemp. Math., 147, Amer. Math. Soc., Providence, RI, 1993.


\end{thebibliography}

\def\JCTB{{\it J.~Combin.\ Theory Ser.\ B}}
\def\CMUC{{\it Comment. Math. Univ. Carol.}}
\def\TAMS{{\it Trans.\ Amer.\ Math.\ Soc.}}
\def\JAMS{{\it J.~Amer.\ Math.\ Soc.}}
\def\PAMS{{\it Proc. Amer. Math. Soc.}}
\def\DM{{\it Discrete Math.}}
\def\CM{{\it Contemporary Math.}}
\def\GC{{\it Graphs and Combin.}}
\def\COM{{\it Combinatorica}}
\def\JGT{{\it J.~Graph Theory}}
\def\JAlgorithms{{\it J.~Algorithms}}
\def\SIAMDM{{\it SIAM J.~Disc.\ Math.}}
\def\CPC{{\it Combinatorics, Probability and Computing}}
\def\EJC{Electron.\ J.~Combin.}
\def\EuropJC{Europ.\ J.~Combin.}

\end{document}